\documentclass[10pt]{article}
\usepackage{amsmath, amsthm}
\usepackage{graphicx}
\newtheoremstyle{theorem}
{10pt} 
{10pt} 
{\sl} 
{\parindent} 
{\bf} 
{. } 
{ } 
{} 
\theoremstyle{theorem}
\newtheorem{theorem}{Theorem}

\newtheorem{lemma}[theorem]{Lemma}
\newtheorem{proposition}[theorem]{Proposition}

\newtheoremstyle{defi}
{10pt} 
{10pt} 
{\rm} 
{\parindent} 
{\bf} 
{. } 
{ } 
{} 
\theoremstyle{defi}

\begin{document}

\title{Characterization of  Cobweb Posets as KoDAGs}
\author{Ewa Krot-Sieniawska \\
Institute of Computer Science,  Bia{\l}ystok University   \\
PL-15-887 Bia{\l}ystok, st.Sosnowa 64, POLAND\\
e-mail: ewakrot@wp.pl, ewakrot@ii.uwb.edu.pl}

\maketitle

\begin{abstract}
The characterization of the large family of cobweb posets  as DAGs
and oDAGs is given. The {\em dim} 2 poset such that its Hasse
diagram coincide with digraf of  arbitrary cobweb poset $\Pi$ is constructed.\\

\vspace{2mm}

\end{abstract}
\vspace{2mm}
\noindent  Presented at Gian-Carlo Rota Polish Seminar:\\
\noindent \emph{http://ii.uwb.edu.pl/akk/sem/sem\_rota.htm}

\section{Introduction}
The   family of the so called cobweb posets $\Pi$ has been
invented by A.K.Kwa\'sniewski in \cite{44,46}. These structures
are such a generalization of the Fibonacci tree growth  that
allows joint combinatorial interpretation for all of them under
the admissibility condition (for references to Kwa\'sniewski
papers see recent \cite{49,49a}).

\noindent Let $\{F_n\}_{n\geq 0}$ be a natural numbers valued sequence with
$F_0=1$ (with $F_0=0$ being exceptional as in case of Fibonacci
numbers). Any sequence satisfying this property uniquely
designates cobweb poset  defined as follows.

\noindent For $s\in\bf{N}_0=\bf{N}\cup\{0\}$ let
$$\Phi_{s}=\left\{\langle j,s \rangle ,\;\;1\leq j \leq F_{s}\right\},\;\;\;$$
Then corresponding cobweb poset is  an infinite partially ordered set
$\Pi=(V,\leq)$, where
 $$ V=\bigcup_{0\leq s}\Phi_s$$
are the elements ( vertices) of $\Pi$ and the partial order relation $\leq$ on $V$ for
 $x=\langle s,t\rangle, y=\langle u,v\rangle $ being  elements of
cobweb poset $\Pi$ is defined by  formula
$$ ( x \leq_{P} y) \Longleftrightarrow
 [(t<v)\vee (t=v \wedge s=u)].$$

\noindent Obviously as any poset a cobweb poset can be represented, via its Hasse
diagram and here it is an infinite directed  graf  $\Pi=\left( V,E\right)$,
\noindent where  set $V$ of its vertices is defined as above and where

$$E =\{\left(\langle j , p\rangle,\langle q ,(p+1) \rangle
\right)\}\;\cup\;\{\left(\langle 1 , 0\rangle ,\langle 1 ,1
\rangle \right)\},$$

\noindent and  where $1 \leq j \leq {F_p}$ and $1\leq
q \leq {F_{(p+1)}}$ stays for  set of (directed) edges.

\noindent For example the Hasse diagram of Fibonacci cobweb poset designated
by the famous Fibonacci sequence looks as follows.

\begin{center}
\includegraphics[width=100mm]{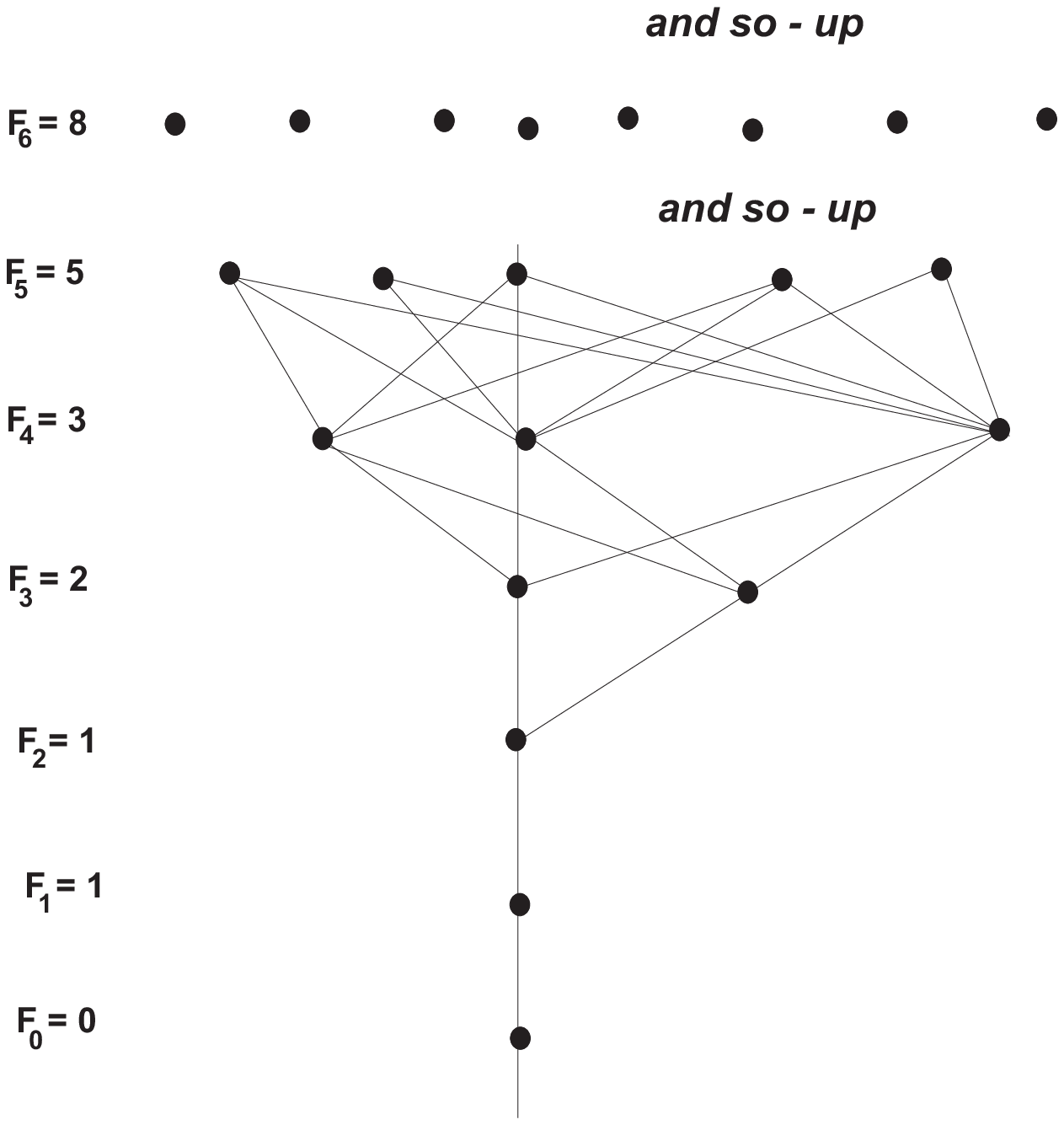}

\vspace{2mm}

\noindent {\small Fig.~1. The construction of the  Fibonacci
"cobweb"  poset}
\end{center}

\noindent The Kwa\'sniewski cobweb posets under consideration
represented by  graphs  are examples of oderable directed acyclic
graphs (oDAG)  which we start to call from now in brief:  KoDAGs.
These are  structures of universal importance for the whole of
mathematics - in particular for discrete "`mathemagics"'
[http://ii.uwb.edu.pl/akk/ ]  and computer sciences in general
(quotation from \cite{49,49a} ):

\begin{quote}
DAGs  considered  as a generalization of trees  have a lot of
applications in computer science, bioinformatics,  physics and
many natural activities of humanity and nature.  For example in
information categorization systems, such as folders in a computer
or in Serializability Theory of Transaction Processing Systems and
many others. Here we introduce specific DAGs as generalization of
trees being inspired by algorithm of the Fibonacci tree growth.
For any given natural numbers valued sequence the graded (layered)
cobweb posets` DAGs  are equivalently representations of a chain
of binary relations. Every relation of the cobweb poset chain is
biunivocally represented by the uniquely designated
\textbf{complete} bipartite digraph-a digraph which is a
di-biclique  designated  by the very  given sequence. The cobweb
poset is then to be identified with a chain of di-bicliques i.e.
by definition - a chain of complete bipartite one direction
digraphs.   Any chain of relations is therefore obtainable from
the cobweb poset chainof complete relations  via
deleting  arcs (arrows) in di-bicliques.\\
Let us underline it again : \textit{any chain of relations is
obtainable from the cobweb poset chain of  complete relations via
deleting arcs in di-bicliques of the complete relations chain.}
For that to see note that any relation  $R_k$ as a subset of  $A_k
\times A_{k+1}$ is represented by a  one-direction bipartite
digraph  $D_k$.  A "complete relation"  $C_k$ by definition is
identified with its one direction di-biclique graph $d-B_k$.  Any
$R_k$ is a subset of  $C_k$. Correspondingly one direction digraph
$D_k$ is a subgraph of an one direction digraph of $d-B_k$.\\
The one direction digraph of  $d-B_k$ is called since now on
\textbf{the di-biclique }i.e. by definition - a complete bipartite
one direction digraph.   Another words: cobweb poset defining
di-bicliques are links of a complete relations' chain.
\end{quote}

\section{DAG $\longrightarrow$ oDAG problem }
In \cite{p} Anatoly D. Plotnikov considered the so called "DAG
$\longrightarrow$ oDAG problem". He had determined condition when a
digraph $G$ may may be presented by the corresponding {\em dim } 2
poset $R$ and he had established the algorithm for how to find it.

\noindent Before citing Plotnikov's results let us however recall
some indispensable definitions following thesource reference \cite{p}.

\noindent Let $P$ and $Q$ be some  partial orders on the same set $A$. Then $Q$
is said to be an {\bf extension} of $P$ if $a\leq_{P} b$ implies
$a\leq_{Q} b$, for all $a,b\in A$. A poset $L$ is a {\bf chain},
or a {\bf linear order} if we have either $a\leq_{L} b$ or
$b\leq_{L} a$ for any $a,b\in A$. If $Q$ is a linear order then it
is a {\bf linear extension} of $P$.

\noindent The {\bf dimension} $dim\ R$ of $R$ being a partial order is the
least positive integer $s$ for which there exists a family $F=(L_1
,L_2 ,\ldots,L_s)$ of linear extensions of $R$ such that $R=
\bigcap_{i=1}^{s} L_{i}$. A family $F=(L_1,L_2,\ldots,L_s)$ of
linear orders on $A$ is called a {\bf realizer} of $R$ on $A$ if
\[
R=\bigcap_{i=1}^{s} L_{i}.
\]

\noindent We shall use $D_{n}$ to denote the set of all acyclic directed $n$-vertex
graphs without loops and multiple edges. Each digraph ${\vec
G}=(V,{\vec E})\in D_{n}$ will be called {\bf DAG}.

\noindent A digraph ${\vec G}\in D_{n}$ shall be called {\bf orderable
(oDAG)} if there exists are $dim\ 2$ poset such that its Hasse
diagram coincides with the digraph ${\vec G}$.

\noindent Let  ${\vec G}\in D_{n}$ be a digraph, which does not contain the
arc $(v_{i},v_{j})$ if there exists the directed path
$p(v_{i},v_{j})$ from the vertex $v_{i}$ into the vertex $v_{j}$
for any $v_{i}$, $v_{j}\in V$. Such digraph is called {\bf
regular}. Let $D\subset D_{n}$ is the set of all regular graphs.

\noindent Let there be given a regular digraph ${\vec G}=(V,E)\in D$, and let
the chain ${\vec X}$ has three elements $x_{i_{1}}$, $x_{i_{2}}$,
$x_{i_{3}}\in X$ such that $i_{1}<i_{2}<i_{3}$, while simultaneously
there are not paths $p(v_{i_{1}},v_{i_{2}})$ ,$p(v_{i_{2}},v_{i_{3}})$
in the digraph ${\vec G}$ but there exists a path
$p(v_{i_{1}},v_{i_{3}})$. Such representation of graph vertices by
elements of the chain ${\vec X}$ is called the representation in
{\bf inadmissible form}. Another words- the chain ${\vec X}$ Constitutes
the graph vertices in {\bf admissible form}.

\noindent Anatoly Plotnikov then shows that:

\begin{lemma}{\em \cite{p}}\label{l1}
\label{l1} A digraph ${\vec G}\in D_{n}$ may be represented by a
$dim\ 2$ poset if:
\renewcommand{\labelenumi}{(\arabic{enumi})}
\begin{enumerate}
\item there exist two chains ${\vec X}$ and ${\vec Y}$, each of
which is a linear extension of ${\vec G}_{t}$; \item the chain
${\vec Y}$ is a modification of ${\vec X}$ with inversions, which
remove the ordered pairs of ${\vec X}$ that there do not exist in
${\vec G}$.
\end{enumerate}
\end{lemma}
\noindent The above lemma results in the algorithm for finding {\em dim} 2
representation of a given DAG (i.e. corresponding oDAG) while the
following theorem establishes the conditions for constructing it.
\begin{theorem}{\em \cite{p}}\label{t1}
\label{t1} A digraph ${\vec G}=(V,{\vec E})\in D_{n}$ can be
represented by $dim\ 2$ poset iff it is regular and its vertices
can be presented by the chain ${\vec X}$ in admissible form.
\end{theorem}

\section{Kwa\'sniewski Cobweb Posets as KoDAGs}
In this section we show that Kwa\'sniewski cobweb posets are
orderable Directed Acyclic Graphs (oDAGs) hence: KoDAGs.

\noindent Obviously,  arbitrary cobweb poset $\Pi=(V, E)$ defined as above
is a DAG (it is directed acyclic graph without loops and multiple
edges). One can also verify that

\noindent \begin{proposition}
 $\Pi=(V, E)$ is a regular digraph.
 \end{proposition}
\noindent \begin{proof}For two elements $\langle i, n\rangle , \langle j,m\rangle \in V$ a
directed path $p(\langle i, n\rangle , \langle j,m\rangle)\notin
E$ will exist iff $n<m+1$ but then  $(\langle i, n\rangle ,
\langle j,m\rangle)\notin E$ i.e. $\Pi$ does not contain the edge
$(\langle i, n\rangle , \langle j,m\rangle)$.
\end{proof}

\noindent It is also possible to verify that vertices of cobweb poset $\Pi$
can be presented in admissible form by the chain ${\vec X}$ being
a linear extension of cobweb $P$ as follows:
\begin{multline*}{\vec X}=\Big(\langle 1,0\rangle,\langle 1,1\rangle ,
\langle 1,2\rangle, \langle 1,3\rangle, \langle 2,3\rangle,
\langle 1,4\rangle, \langle 2,4\rangle, \langle 3,4\rangle,\langle
1,5\rangle, \langle 2,5\rangle, \langle 3,5\rangle,\\\langle
4,5\rangle,\langle 5,5\rangle,...\Big),\end{multline*}
 where
$$ ( \langle s,t\rangle \leq_{{\vec X}} \langle u,v\rangle) \Longleftrightarrow
 [(t < v)\vee (t=v \wedge s\leq u)]$$

\noindent for $1\leq s \leq F_{t},\; 1\leq u \leq F_{v},\;\;\;t, v \in {\bf N}\cup\{0\}.$

\noindent Then the cobweb poset $\Pi$ satisfies conditions of the Theorem
\ref{t1}, so it is oDAG. In order to find the chain ${\vec Y}$ being
a linear extension of cobweb $P$ one uses the Lemma \ref{l1} thus
arriving  at: \begin{multline*} {\vec Y}=\Big(\langle
1,0\rangle,\langle 1,1\rangle , \langle 1,2\rangle, \langle
2,3\rangle, \langle 1,3\rangle, \langle 3,4\rangle, \langle
2,4\rangle, \langle 1,4\rangle,\langle 5,5\rangle, \langle
4,5\rangle, \langle 3,5\rangle,\\\langle 2,5\rangle,\langle
1,5\rangle,...\Big),\end{multline*}
\noindent where

$$ ( \langle s,t\rangle \leq_{{\vec Y}} \langle u,v\rangle) \Longleftrightarrow
 [(t < v)\vee (t=v \wedge s\geq u)]$$
 for $1\leq s \leq F_{t},\; 1\leq u \leq F_{v},\;\;\;t, v \in {\bf
 N}\cup\{0\}$   and finally
 $$ (P,\leq_{P})={\vec X}\cap{\vec Y}.$$

\noindent So we have the announced charcterization Theorem $4$ being proved now.
\noindent \begin{theorem}
 An arbitrary cobweb poset $\Pi=(V, E)$ is an example od orderable
 directed acyclic digraph (oDAG) and $ \Pi=(V,\leq)={ X}\cap{ Y},$
 for $X$, $Y$ being linear extensions of partial order $\Pi$ defined as above.
\end{theorem}

\vspace{2mm}

\noindent \textbf{Ackowledgements. } My  thanks are to  Professor
A. Krzysztof Kwa\'sniewski for guidance and final improvements of
this paper. Attendance of W. Bajguz and Maciej Dziemianczuk is
highly appreciated too.

\end{document}